\documentclass[12pt]{article}
\usepackage{ amsmath, amsthm, amssymb, amsfonts}
\usepackage[T1]{fontenc}
\usepackage[latin1]{inputenc}
\usepackage[extra]{tipa} 
\newcommand\ov{\over}

\renewcommand{\s}{\sigma} 
\renewcommand\t{\tau}

\newcommand\ve{\varepsilon}
\newcommand{\la}{\lambda}

\newcommand{\x}{\xi}

 \def\({\left(} \def\){\right)}

\newcommand{\ra}{\rightarrow}
 \newcommand{\bE}{\mathbb{E}}
\newcommand{\bP}{\mathbb{P}}

\newcommand{\bZ}{\mathbb{Z}}

\newcommand{\cC}{\mathcal{C}}

\newcommand{\cI}{\mathcal{I}}
\newcommand{\cJ}{\mathcal{J}}

\newcommand{\be}{\begin{equation}}
\newcommand{\ee}{\end{equation}}
\newcommand{\ba}{\begin{eqnarray*}}
\newcommand{\ea}{\end{eqnarray*}}
\newcommand{\bae}{\begin{eqnarray}}
\newcommand{\eae}{\end{eqnarray}}
\newcommand{\bc}{\begin{center}}
\newcommand{\ec}{\end{center}}

\newcommand{\qBi}[2]{\bigg[{#1\atop #2}\bigg]}
\newcommand{\fr}{\frac}

\begin{document}
\begin{titlepage}
\hfill July 20, 2009
\begin{center}{\Large\bf  On the Distribution of a Second Class Particle\\ in the
\vspace{.75ex}
 \\Asymmetric Simple Exclusion Process}\end{center}
 
\begin{center}{\large\bf Craig A.~Tracy}\\
{\it Department of Mathematics \\
University of California\\
Davis, CA 95616, USA\\
email: \texttt{tracy@math.ucdavis.edu}}\end{center}

\begin{center}{\large \bf Harold Widom}\\
{\it Department of Mathematics\\
University of California\\
Santa Cruz, CA 95064, USA\\
email: \texttt{widom@ucsc.edu}}\end{center}

\begin{abstract}
We give an exact expression for the distribution of the position $X(t)$ of a single
 second class particle in the asymmetric simple exclusion
process (ASEP) where initially the second class particle is located  at the origin and the first class particles occupy
the sites $\bZ^+=\{1,2,\ldots\}$.
\end{abstract}
\end{titlepage}
 \section{Introduction}
The asymmetric simple exclusion process (ASEP)  \cite{Li1, Li2} is one of the simplest models of nonequilibrium statistical
mechanics  and has been called the ``default stochastic model for transport phenomena'' \cite{Yau}.  
A useful concept
in exclusion processes is that of a \textit{second class particle}:\footnote{The following quote is taken from Liggett \cite{Li2}.}

\begin{quote}
Imagine that the particles in the system are each called either first class or
second class.  The evolution is the same as before, except that if a second
class particle attempts to go to a site occupied by a first class particle, it
is not allowed to do so, while if a first class particle attempts to move to a site
occupied by a second class particle, the two particles exchange positions.
In other words, a first class particle has priority over a second class particle.
This rule has no effect on whether or not a give site is occupied at a given time.
The advantage, though, is that viewed by itself, the collection of first class particles
is Markovian, and has the same law as the exclusion process.  The collection
of second class particles is clearly not Markovian.  However, the collection
of first and second class particles is Markovian, and again evolves like an
exclusion process.
\end{quote}
Here we consider ASEP on the integer lattice $\bZ$ with jumps
one step to the right with rate $p$ and jumps one step to the left with rate $q=1-p$.  We assume a leftwards drift, i.e.\ $q>p$.  We 
further assume that the system has one second class particle initially located at the origin and first class particles initially located at sites in
\[ Y=\{0<y_1<y_2<\cdots \}\subset \bZ^+.\]
With the above initial condition, we denote by $X(t)$  the position of the second class particle at time $t$.  The purpose of this note is to give an exact expression for
the probability that the second class particle is at position $x$ at time $t$, i.e.\
$\bP_Y\left(X(t)=x\right)$.  (The subscript $Y$ denotes the sites of the initial configuration
of the first class particles.)  Our main result is for $Y=\bZ^+$ and is given below in (\ref{J3sum}) and in
a slightly different form in (\ref{distrFn}).

\section{A Basic Lemma}
The single second class particle located at $X(t)$  can be viewed as the (single) discrepancy under \textit{basic coupling} between two  asymmetric simple exclusion processes $\eta_t$ and $\zeta_t$  where 
$\zeta_t(X(t))=1$ and $\eta_t(X(t))=0$ and initially 
$\left\{x:\zeta_0(x)=1\right\}=Y'=\{0\}\cup Y$ and 
$\left\{x:\eta_0(x)=1\right\}=Y$ \cite{Li1, Li2}. 

The following identity we first learned from H.~Spohn \cite{Spo} but presumably
it has a long history:
\be \bP_Y\left(X(t)=x\right)=\bP_{Y'}\left(\zeta_t(x)=1\right)-\bP_Y\left(\eta_t(x)=1\right).
\label{lemma}\ee

For the convenience of the reader, we give a short proof of (\ref{lemma}).  Let $\zeta_t$
and $\eta_t$  be as above evolving together under
the basic coupling \cite{Li1, Li2}. Recall that the coupled processes
satisfy $\eta_t\le \zeta_t$ for all $t>0$ since they satisfy this inequality at $t=0$
 \cite{Li1, Li2}.\footnote{Given two configurations $\eta,\zeta\in\{0,1\}^{\bZ}$ we say $\eta\le\zeta$ if $\eta(x)\le \zeta(x)$ for all $x\in\bZ$.}
Define
\ba
\cJ_\eta(x,t)&:=& \sum_{z\le x} \eta_t(z)=\textrm{number of particles in
configuration}\>\eta_t\>\>\textrm{with positions}\>\le x, \\
\cJ_\zeta(x,t)&:=& \sum_{z\le x} \zeta_t(z)
=\textrm{number of particles in
configuration}\>\zeta_t\>\>\textrm{with positions}\>\le x, \\
\cI(x,t)&=&\left\{\begin{array}{ll} 1 & \textrm{if}\>\> X(t)\le x,\\
						0 & \textrm{if}\>\> X(t)>x.\end{array}\right.\ea
By counting
\be \cJ_{\zeta}(x,t)= \cJ_\eta(x,t)+\cI(x,t). \label{count}\ee
Since
\[ \bE_{Y'}\left(\cJ_{\zeta}(x,t)\right)=\sum_{z\le x} \bE_{Y'}\left(\zeta_t(z)\right)
=\sum_{z\le x}
\bP_{Y'}(\zeta_t(z)=1),\]
\[ \bE_Y\left(\cJ_{\eta}(x,t)\right)=\sum_{z\le x} \bE_Y\left(\eta_t(z)\right)=\sum_{z\le x}
\bP_Y(\eta_t(z)=1),\]
\[ \bE_Y\left(\cI(x,t)\right)=\bP_Y\left(X(t)\le x\right)=\sum_{z\le x} \bP_Y(X(t)=z),\]
the expectation of (\ref{count}) gives
\[ \sum_{z\le x} \bP(X(t)=z)=\sum_{z\le x}
\bP(\zeta_t(z)=1)-\sum_{z\le x}
\bP(\eta_t(z)=1)\]
from which (\ref{lemma}) follows.
\section{Probability for a site to be occupied in ASEP}
For ASEP with particles initially at $Y$ we denote by $x_m(t)$ the position
of the $m$th left-most particle at time $t$ (so $x_m(0)=y_m$).  In
Theorem 5.2 of  \cite{TW1} the authors gave an exact expression for 
$\bP_Y(x_m(t)=x)$.
To state this result we first
 recall the definition of the $\t$-binomial coefficients. For $0\le \t:=p/q <1$ define for each $n\in\bZ^+$
\[ [n]={1-\t^n\ov 1-\t},\>\>\> [n]!=[n][n-1]\cdots [1],\> [0]!:=1,\>\>\> 
\qBi{n}{k}={[n]!\ov [k]! [n-k]!}, \>0\le k \le n,\]
and if $k>n$ we set $\left[{n\atop k}\right]=0$.  Equation (5.12) of \cite{TW1} can be
written in the following way\footnote{We make some changes in the notation in (5.12) of
\cite{TW1}.  The $(p,q)$-binomial coefficient $\left[{n\atop k}\right]$ of \cite{TW1} equals $q^{k(n-k)}$ times the 
$\t$-binomial coefficient $\left[{n\atop k}\right]$ defined above. The second change
is a little more subtle.  The sum in (5.12) is over all finite subsets 
$S\subset\{1,2,\ldots,|Y|\}$
with $|S|\ge m$.
If $S=\{s_1,\ldots, s_k\}$ the subset $Y_S:=\{y_{s_1},\ldots, y_{s_k}\}$
and the factor $\prod_{i\in S} \x_i^{-y_i}$ appears in the integrand of (5.12). Thus
we can equivalently sum over all finite subsets $S\subset Y$ where now the factor
$\prod_{1\le i\le k} \x_i^{-s_i}$ appears in the integrand. The factor
$\s(S)=\sum_{i\in S} i$ of (5.12) becomes $\s(S,Y)$ given above.}$^{,}$\footnote{All contour integrals are to be given a factor of $1/2\pi i$.}
\be
\bP_Y\left(x_m(t)=x\right)=\sum_{k=1}^{|Y|}
\sum_{{S\subset Y\atop|S|=k}} c_{m,k}\, \t^{\s(S,Y)}\int_{\cC_R}
\cdots \int_{\cC_R} I(x,k,\x)\, \prod_{i=1}^k \x_i^{-s_i}\, d^{k}\x \label{position}\ee
where,  if $S:=\{s_1,\ldots, s_k\}$ then
\ba
c_{m,k}&=& q^{k(k-1)/2} (-1)^{m+1} \t^{m(m-1)/2} \t^{-km} \qBi{k-1}{k-m}, \\
\s(S,Y)&=&\#\left\{(s,y): s\in S, y\in Y,\>\textrm{and}\> y\le s\right\}\\
&=& \textrm{sum of the positions of the elements of}\>\> S\>\textrm{in}\>\> Y,\\
I(x,k,\x)&=&\prod_{1\le i<j\le k}{\x_j-\x_i\ov p+q\x_i\x_j-\x_i}\,
\left(1-\prod_{i=1}^k \x_i\right) \prod_{i=1}^k
{\x_i^{x-1} e^{\ve(\x_i) t}\ov 1-\x_i}\,, \\
\ve(\x)&=&\fr{p}{\x}+q\,\x -1 
\ea
and $\cC_R$ is a circle of radius $R$ centered at the origin with $R\gg 1$ so that
all (finite) singularities of the integrand are enclosed by $\cC_R$.  Observe that
$c_{m,k}=0$ when $m>k$.

Since
\be \bP_Y\left(\eta_t(x)=1\right)=\sum_{m=1}^{|Y|} \bP_Y\left(x_m(t)=x\right),
\label{density}\ee
we sum the right side of (\ref{position}) over all $m\le k$. To carry out this sum
recall the $\t$-binomial theorem
\[\sum_{j=0}^n \qBi{n}{j} (-1)^j z^j \t^{j(j-1)/2}=(1-z)(1-z\t)\cdots (1-z\t^{n-1}).\]
Using this a simple calculation shows
\[ \sum_{m=1}^k (-1)^{m+1} \t^{m(m-1)/2} \t^{-km} \qBi{k-1}{k-m}=(-1)^{k+1}
 \t^{-k(k+1)/2} \prod_{j=1}^{k-1} (1-\t^j). \]
 Thus
 \bae \bP_Y\left(\eta_t(x)=1\right)&=&\sum_{k=1}^{|Y|} (-1)^{k+1} q^{k(k-1)/2}
 \t^{-k(k+1)/2} \prod_{j=1}^{k-1}(1-\t^j)\nonumber \\
 && \quad \times 
 \sum_{{S\subset Y\atop|S|=k}} \t^{\s(S,Y)}\int_{\cC_R}
\cdots \int_{\cC_R} I(x,k,\x)\, \prod_{i=1}^k \x_i^{-s_i}\, d^{k}\x. \label{Isum}\eae

\noindent{\textsc{Remark:}} The above formula holds for
either $|Y|$ finite or infinite.  For $|Y|=N$,
the integral of order $N$ in (\ref{Isum}) is gotten from the summand $S=Y$.
Since $\s(Y,Y)=N(N+1)/2$, we get for the coefficient of this integral
\be (-1)^{N+1} q^{N(N-1)/2} \prod_{j=1}^{N-1}=(-1)^{N+1} \prod_{j=1}^{N-1}(q^j-p^j).
\label{coeff}\ee
\section{Probability for a site to be occupied by a second class particle}
As above, suppose that our initial configuration consists of a second class particle at site 0 and first class particles at sites in $Y$. As above, set $Y'=Y\cup \{0\}$. 
The process $\zeta_t$ has initially its particles at sites in $Y'$.  
 We apply formula (\ref{Isum}) to the initial configurations $Y'$ and $Y$ and by
 (\ref{lemma}) we subtract to obtain $\bP_Y(X(t)=x)$.
  If $|Y'|=N$ there is one $N$-dimensional integral that comes from the expansion of  $\bP_{Y'}(\zeta_t(x)=1)$ when $S=Y'$. The coefficient of the integral of highest order equals (\ref{coeff}). 
  
  We now consider the special case of \textit{step initial condition}; that is,
  $Y=\bZ^+$, and use Corollary (5.13) of \cite{TW1} to obtain a more compact
  expression for $\bP_{\bZ^+}(x_m(t)=x)$.  To find $\bP_{\bZ^+}(\eta_t(x)=1)$
  we again apply (\ref{density}) but use (5.13) of \cite{TW1}.  As above we interchange
  the sums over $k$ and $m$, use the $\t$-binomial theorem (\cite{Mac}, pg.~26), to conclude
  \be
 \bP_{\bZ^+}\hspace{-.5ex}\left(\eta_t(x)=1\right)=-\sum_{k\ge 1} \fr{q^{k^2}}{k!}
 \prod_{j=1}^{k-1}(1-\t^j)
 \int_{\cC_R}\cdots \int_{\cC_R} \tilde{J}_k(x,\x)\,d\x_1\cdots d\x_k
 \label{Jsum} \ee
 where
 \[\tilde{J}_k(x,\x)=\prod_{i\ne j}{\x_j-\x_i\ov p+q\x_i\x_j-\x_i}\;(1-\prod_i\x_i)\;\prod_i{\x_i^{x-1}\,e^{\ve(\x_i)t}\ov (1-\x_i)\,(q\x_i-p)}.\]
We can get  the corresponding formula for $Y'=\bZ^+\cup \{0\}$ by observing that there is a one-one correspondence between subsets $S'\subset Y'$ and subsets $S\subset Y$ given by $S=S'+1$. Then 
$\s(S',\,Y')=\s(S,\,Y)$ and, with obvious notation, $\prod \x_i^{-{s_i}'}=\prod\x_i\,\cdot\,\prod \x_i^{-s_i}$. It follows that for the difference 
$\bP_{Y'}(\zeta_t(x)=1)-\bP_Y(\eta_t(x)=1)$ we multiply the integrand $\tilde{J}_k(x,\x)$ in (\ref{Jsum}) by $\prod\x_i-1$. 

Thus
\be \bP_{\bZ^+}\hspace{-.5ex}\left(X(t)=x\right)=\sum_{k\ge 1} \fr{q^{k^2}}{k!}
 \prod_{j=1}^{k-1}(1-\t^j)
 \int_{\cC_R}\cdots \int_{\cC_R} \doubletilde{J}_k(x,\x)\,d\x_1\cdots d\x_k
 \label{J2sum} \ee
 where
 \[\doubletilde{J}_k(x,\x)=\prod_{i\ne j}{\x_j-\x_i\ov p+q\x_i\x_j-\x_i}\:(1-\prod_i\x_i)^2\:\prod_i{\x_i^{x-1}\,e^{\ve(\x_i)t}\ov (1-\x_i)\,(q\x_i-p)}.\]
 From this it follows that the distribution function is (on $\cC_R$,  $|\x^{-1}|\ll 1$)
 \be
  \bP_{\bZ^+}\hspace{-.5ex}\left(X(t)\le x\right)=\sum_{k\ge 1} \fr{q^{k^2}}{k!}
 \prod_{j=1}^{k-1}(1-\t^j)
 \int_{\cC_R}\cdots \int_{\cC_R} J_k(x,\x)\,d\x_1\cdots d\x_k
 \label{J3sum} \ee
  where
 \[J_k(x,\x)=\prod_{i\ne j}{\x_j-\x_i\ov p+q\x_i\x_j-\x_i}\;(\prod_i\x_i-1)\,\prod_i{\x_i^{x}\,e^{\ve(\x_i)t}\ov (1-\x_i)\,(q\x_i-p)}.\]
Since
\[ {1\ov p + q\x\x' -\x}={1\ov \x (\x'-1)}+\textrm{O}(\t),\>\>\t\ra 0,\]
the TASEP limit of $J_k(x,\x)$ is
\[ J_k^{\textrm{\tiny{TASEP}}}(x,\x):=\lim_{\t\ra 0}J_k(x,\x)=\prod_{i\neq j}(\x_j-\x_i)\, (\prod \x_i -1) \prod_i {\x_i^x e^{\ve(\x_i)t}\ov
\left(\x_i(1-\x_i)\right)^k} \]
where now $\ve(\x)=\x-1$; and hence,
\be \lim_{\t\ra 0}\bP_{\bZ^+}\hspace{-.5ex}\left(X(t)\le x\right)=\sum_{k\ge 1} \fr{1}{k!}
\int_{\cC_R}\cdots \int_{\cC_R} J^{\textrm{\tiny{TASEP}}}_k(x,\x)\,d\x_1\cdots d\x_k.
 \label{TASEPdistr} \ee
Expression (\ref{J3sum}) for the distribution function
 can be simplified somewhat.  Define the kernel
\[ K_{x,t}(\x,\x')=q\,{(\x')^x e^{\ve(\x')t}\ov p + q\x\x'-\x},\]
and the associated operator $K_{x,t}$ on $L^2(\cC_R)$ by
\[ f(\x)\longrightarrow \int_{\cC_R} K_{x,t}(\x,\x') f(\x')\, d\x',\> \x\in\cC_R. \]
Then using the identity \cite{TW2}
\[ \det\left({1\ov p +q\x_i\x_j-\x_i}\right)_{1\le i,j\le k}
= (-1)^k(pq)^{k(k-1)/2} \prod_{i\ne j} {\x_j-\x_i\ov p +q \x_i\x_j -\x_i}\, \prod_i{1\ov (1-\x_i)(q\x_i-p)}\]
we have
\bae\hspace{-3ex} \bP_{\bZ^+}\hspace{-.2ex}(X(t)\le x)
&\hspace{-1ex}=&\sum_{k\ge 1} \t^{-k(k-1)/2} \prod_{j=1}^{k-1}(1-\t^j)
\times\nonumber\\
&& \fr{(-1)^k}{k!}\int_{\cC_R}\cdots \int_{\cC_R} \left[ \det\left(K_{x+1,t}(\x_i,\x_j)\right)_{1\le i,j\le k}-\det\left(K_{x,t}(\x_i,\x_j)\right)_{1\le i,j\le k}\right]\nonumber\\
&\hspace{-10ex}=&\hspace{-5ex}\sum_{k\ge 1} \t^{-k(k-1)/2}\prod_{j=1}^{k-1}(1-\t^j)
\int_{\cC_R}\fr{1}{\la^{k+1}}\left[\det(I-\la K_{x+1,t})-\det(I-\la K_{x,t})\right]\, d\la
\label{distrFn}\eae
where $\det(I-\la K_{x,t})$ is the Fredholm determinant and the last line follows from the Fredholm expansion.

\pagebreak[2]
\noindent{\textsc{Remarks:}}
\vspace{-2.5ex}
\begin{enumerate}
\item One cannot interchange the sum and the integration in (\ref{distrFn})  as was possible in an analogous calculation in \cite{TW2}. 
This is the case even though (\ref{J3sum}) converges absolutely  for
all $0\le \t \le 1$ (recall one may take $R\gg 1$).  Thus we do not have a representation of
$\bP_{\bZ^+}(X(t)\le x)$ as a single integral whose integrand involves the above
 Fredholm determinants as was the case in \cite{TW2}. 
\item ASEP with first and second class particles is integrable in the sense that the Yang-Baxter equations are satisfied \cite{AB}.  Using this integrable structure, it is possible
to compute directly, i.e.\ without using the basic lemma (\ref{lemma}), 
$\bP_{\bZ^+}(X(t)=x)$  using methods similar to that of
 \cite{TW1}.  We have carried this out to the extent that (\ref{coeff}) was computed by this approach. However, this route is much more involved than the one presented here.
\end{enumerate}

\noindent{\textsc{Acknowledgements:}} The authors wish to thank Professor Dr.~Herbert Spohn for communicating to them the basic lemma.  The authors thank Professor Dr.~Gert-Martin Greul and the staff at the Mathematisches Forschungsinstitut Oberwolfach for their
hospitality during the authors' Research in Pairs stay.  This work was supported
by the National Science Foundation under grants DMS--0553379 (first author)
and DMS--0552388 (second author).

\end{document}